\documentclass{amsart}
\usepackage{doc}
\usepackage{hyperref}
\usepackage{amsmath, amsthm, amssymb, amsfonts}
\usepackage{nicefrac}
\newtheorem{theorem}{Theorem}[section]
\newtheorem{lemma}[theorem]{Lemma}
\newtheorem{definition}[theorem]{Definition}
\newtheorem{remark}[theorem]{Remark}
\newcommand{\N}{\mathbb{N}}
\newcommand{\Z}{\mathbb{Z}}
\newcommand{\B}{\mathrm{B}}
\newcommand{\R}{\mathbb{R}}
\newcommand{\C}{\mathbb{C}}
\newcommand{\E}[1]{\mathbb{E}\left[#1\right]}
\newcommand{\rL}{\right\}}
\newcommand{\lL}{\left\{}
\newcommand{\rC}{\right ]}
\newcommand{\lC}{\left [}
\newcommand{\rP}{\right)}
\newcommand{\lP}{\left(}

\newcommand{\Prob}[1]{\mathbb{P}\left(#1\right)}
\newcommand{\abs}[1]{\lvert#1\rvert}
\makeatletter
\def\@setthanks{\vspace{-\baselineskip}\def\thanks##1{\@par##1}\thankses}
\makeatother
\begin{document}
\title[Zero-free neighborhood for Kac polynomials]{Zero-free neighborhood around the
Unit Circle for Kac polynomials}
\author{Gerardo Barrera*}
\address{University of Helsinki, Department of Mathematical and Statistical Sciences.
Exactum in Kumpula Campus, Pietari Kalmin katu 5.
Postal Code: 00560. Helsinki, Finland.}
\email{gerardo.barreravargas@helsinki.fi}
\thanks{*Corresponding author: Gerardo Barrera. Email: gerardo.barreravargas@helsinki.fi}
\thanks{G. Barrera was supported by CIMAT and PIMS.}
\author{Paulo Manrique}
\address{General Coordination of Institutional Organization and Information.
National Polytechnic Institute.
Postal code: 07738. Mexico city, Mexico.}
\email{pmanriquem@ipn.mx}
\thanks{P. Manrique was supported by C\'atedras-CONACyT, M\'exico.}
\subjclass[2000]{Primary 60G99,  12D10; Secondary, 11CXX, 30C15}
\keywords{Locally Sub-Gaussian random variables; 
Salem-Zygmund type inequalities;
Small ball probability; Zeros of random polynomials}

\begin{abstract}
In this paper, we study how the roots of the so-called Kac polynomial $W_n(z) = \sum_{k=0}^{n-1} \xi_k z^k$ are concentrating to the unit circle when its coefficients of $W_n$ are independent and identically distributed non-degenerate real random variables.  It is well-known that the roots of a Kac polynomial are concentrating around the unit circle as $n\to\infty$ if and only if 
$\mathbb{E}[\log( 1+ |\xi_0|)]<\infty$.  Under the condition of $\mathbb{E}[\xi^2_0]<\infty$, we show that there exists an annulus of width 
$\textnormal{O}(n^{-2}(\log n)^{-3})$ around the unit circle which is {\it free} of roots with probability 
$1-\textnormal{O}({(\log n)^{-{1}/{2}}})$.
The proof relies on the so-called small ball probability inequalities and the least common denominator 
used in \cite{RV2}.
\end{abstract}
\maketitle

\section{Introduction}\label{introduction}
\markboth{G. Barrera and P. Manrique}{Zero free neighborhood for Kac polynomials}

In the signal analysis and the speech recognition, an important tool is the so-called $z$-transform (a particular case is the discrete Fourier transform). 
In this context, the study of zeros of the $z$-transform provides useful information of a signal. The existence of a region free of zeros around the unit circle is an important aspect for a proper work of the $z$-transform. For further details, see Chapter 3 in \cite{BOZK}. 
Roughly speaking, in this paper we find a region free of zeros around the unit circle for the so-called Kac polynomials which are closely related to the discrete Fourier transform.

Let $n\in \mathbb{N}$ and 
let $\xi_0,\ldots,\xi_{n-1}$ be independent and identically distributed (iid for shorthand) non-degenerate real random variables (rvs for shorthand) defined in the probability space $(\Omega,\mathcal{F},\mathbb{P})$. Denote by $\mathbb{E}$ the expectation respect to the measure $\mathbb{P}$. The Kac polynomial $W_n$ is defined as the random polynomial of degree $n-1$ given by 
\[W_n(z)=\sum_{j=0}^{n-1} \xi_j z^j,\quad z\in \mathbb{C}.
\] 

In the sequel, we introduce the basic notation and terminology which will be used throughout this manuscript.
For any $z\in \mathbb{C}$, denote by $|z|$ and $\arg(z)$ the modulus of $z$ and the argument of $z$, respectively. Let $a,b\in \mathbb{R}$ such that $a\leq b$.
Denote by $R_n(a,b)$ the number of roots of $W_{n}$ in the annulus $\lL z\in\C : a\leq \abs{z}\leq b\rL$ 
and for any $\alpha,\beta\in [-\pi,\pi]$ such that $\alpha\leq \beta $ denote by $S_n(\alpha,\beta)$
the number of roots in $\lL z\in\C: \alpha \leq \arg(z)\leq \beta\rL$. 
 
Shparo and Shur \cite{ShpShu1962} proved that under general conditions on the random coefficients (rcs for shorthand), the roots of $W_n$ concentrate around the unit circle with asymptotically uniform distribution in the argument as $n$ goes by. Moreover, Ibragimov and Zaporozhets \cite{IbrZap2013} showed that if the rcs of $W_n$ are non-degenerate and satisfy $\mathbb{E}[\log(1 + |\xi_0|)] < \infty$ if and only if its roots are asymptotically concentrated near the unit circle. Later, Kabluchko and Zaporozhets \cite{KabZap2012} provided a wide description of the localization of the roots for different conditions on the rcs.  
We point out that the localization of the roots of Kac polynomials determine the poor efficiency of some algorithms from the speech recognition and signal processing applications, see \cite{DrugYann2015} for further details. 

Ibragimov and Zaporozhets \cite{IbrZap2013} proved that
\[
\mathbb{P}\Big(\lim\limits_{n\to \infty}\frac{1}{n} R_n\lP 1-\delta, 1+\delta\rP = 1\Big) = 1\quad \textrm{ holds for any } \delta\in(0,1)\]  
if and only if 
$\mathbb{E}[\log(1 + |\xi_0|)] < \infty$. They also proved that for any distribution $\xi_0$ and $\alpha,\beta\in\lP -\pi, \pi\rP$ the following holds:
\[\mathbb{P}\Big(
\lim\limits_{n\to \infty}\frac{1}{n} S_n\lP \alpha,\beta\rP =  \frac{\beta-\alpha}{2\pi}\Big)= 1.\]
Shepp and Vanderbei \cite{SheVan1995} studied the case of iid standard Gaussian coefficients and showed that
\begin{equation}\label{eqn160920191424} 
\lim\limits_{n\to \infty}\frac{1}{n} \mathbb{E}[
R_n (e^{-{\delta}/{n}}, 
e^{{\delta}/{n}})] = \frac{1+e^{-2\delta}}{1-e^{-2\delta}} - \frac{1}{\delta}\quad\textrm{ for any } \delta >0.
\end{equation}
Later, Ibragimov and Zeitouni\cite{IbrZei1997} extended \eqref{eqn160920191424} to the case of iid coefficients which common distribution belongs to the domain of attraction of 
an $\alpha$-stable law,
\begin{equation}\label{eqn190202191014}
\lim\limits_{n\to \infty}\frac{1}{n} \mathbb{E}[
R_n(e^{-{\delta}/{n}}, 
e^{{\delta}/{n}})]= \frac{1+e^{-\alpha \delta}}{1-e^{-\alpha\delta}} - \frac{2} {\alpha\delta}
\quad \textrm{ for any } \delta >0.
\end{equation} 
Note that for any $\delta>0$, as $\alpha \to 0^+$ we have $\frac{1+e^{-\alpha \delta}}{1-e^{-\alpha\delta}} - \frac{2}{\alpha\delta} \rightarrow 0$.
Then \eqref{eqn190202191014} may tend to zero as $n\to\infty$ when $\xi_0$ has a slowly varying tail distribution. In fact, G\"{o}tze and Zaporozhets \cite{GotZap2012} showed that if $\abs{\xi_0}$ has a slowly varying tail distribution, then 
\begin{equation*}\label{eqn200220191248}
\lim\limits_{n\to \infty}
\mathbb{P}\big(R_n(e^{-{\delta}/{n}}, e^{{\delta}/{n}}) = 0\big) =1\quad \textrm{ for any } \delta >0,
\end{equation*}
i.e., the roots of a Kac polynomial with iid rcs with a slowly varying tail distribution hit the unit circle with almost zero probability. 

In the case that $W_{n}$ has iid rcs which common distribution belongs to the domain of attraction of an $\alpha$-stable law, limit \eqref{eqn190202191014} yields that for $\delta>0$, $W_{n}$ has at least one root in the annulus 
${R}_{\delta,n}:=\{ z\in\C : e^{-{\delta}/{n}} \leq \abs{z}\leq e^{{\delta}/{n}}\}$ with positive probability for all large $n$ and
\[\mathbb{P}\big(R_n(e^{-{\delta}/{n}}, e^{{\delta}/{n}}) = n\big) \leq \frac{1+e^{-\alpha \delta}}{1-e^{-\alpha\delta}} - \frac{2}{\alpha\delta} + \textnormal{o}(1).\] 
Therefore, a remarkable question is to determine if there exists an annulus inside of ${R}_{\delta,n}$ such that $W_{n}$ has at least one root or not on it. 
The existence of roots \textit{pretty close} to the unit circle is an important aspect in the analysis of signals. This helps to understand the contribution of the phase information of a signal.
We refer to \cite{DrugYann2015} for further details.

Shepp and Vanderbei \cite{SheVan1995} conjectured that with high probability the nearest root of $W_n$ to the unit circle is at a distance of order $\textnormal{O}(n^{-2})$. Later, Konyagin and Schlag \cite{Kon1999} showed that the last conjecture holds true when the rcs have standard Gaussian or Rademacher (uniform distribution on $\{-1,1\}$) distribution. To be more precise, Konyagin and Schlag proved that there exists a positive constant $C$ such that for any $t>0$ 
\begin{equation}\label{eq:kay}
\limsup_{n\to\infty} 
\mathbb{P}\Big(\min_{\abs{\abs{z}-1}\leq t n^{-2}} \abs{W_n(z)} \leq t n^{-{1}/{2}} \Big) \leq C t.
\end{equation}
They also showed 
\begin{equation}\label{eq:limitzero} \mathbb{P}\left(
\min_{x\in[0,1]}|W_n(x)|\leq n^{-1/2}(\log n)^{-\gamma}\right)=\textnormal{o}(1), \; \textrm{ as } n\to \infty
\end{equation}
for $\gamma>1/2$ and iid Gaussian rcs.

Karapetyan \cite{Kar1998} mentioned that it is possible to extend the above result under the assumption of non-degenerate real {\it{sub-Gaussian}} rcs,
but only a sketch of the proof was given.
Moreover, he claimed that the previous result can be extended under the finiteness of the third moment on the rcs.
However,
Karapetyan \cite{Karapetyan2020} showed that for iid rcs with zero mean
and finite third moment, it follows
for any $\epsilon \in (0,1)$
and $n>16C^\frac{9936}{\epsilon^3}$,
\begin{equation}\label{eq:kay1}
\mathbb{P}\Big(
\min\limits_{x\in [0,1]}
\big|
\sum\limits_{j=0}^{n-1}\xi_j e^{ijx}
\big|
\geq n^{-1/2+\epsilon}\Big)\leq 
\frac{1}{n^{\nicefrac{\epsilon^2}{180}}},
\end{equation}
where the constant $C$ depends only on the moments of  $\xi_0$.
The proof of \eqref{eq:kay1}
is long, technical and complicated.

Later, Barrera and Manrique \cite{GP} proved that if the moment generating function of iid coefficients exists in an open neighborhood around $0$, then for any $t\geq 1$ 
\begin{equation}\label{chin}
\mathbb{P}\Big(\min_{\abs{\abs{z}-1} \leq  tn^{-2}\lP\log n \rP^{-{1}/{2}-\gamma}} \left|W_n(z)\right| \leq t n^{-{1}/{2}}(\log n)^{-\gamma}\Big)= \textnormal{O}((\log n)^{-\gamma+{1}/{2}}),
\end{equation}
where $\gamma>{1}/{2}$. 
The proof of \eqref{chin} recovers the essential ideas of Konyagin and Schlag \cite{Kon1999} who only considered the problem when the rcs have Rademacher or standard Gaussian distribution. 
Their proof is quite technical and involved. It is based on the so-called Salem-Zygmund inequality for sub-Gaussian rvs.

To extend for more distributions, Barrera and Manrique \cite{GP} took advantage of the concept of {\it least common denominator} (lcd for shorthand), which was developed in the study of the singularity of random matrices \cite{RV2}. 
Roughly speaking, the lcd is a combinatorial measurement to understand the concentration of a sum of independent rvs in a small ball.
Furthermore, under the assumptions of the finiteness of the second moment,  using similar ideas from Barrera and Manrique \cite{GP}, it is possible to find an annulus in which $W_n$ does not have roots with high probability. 
We remark that the lcd has been converted into a useful tool that allows analyzing different interesting problems. For instance, it is used in the study of isomorphism between graphs \cite{Luh2018} and in the analysis of the condition number for random matrices \cite{RV1}.
In this manuscript, the lcd is used to understand how small can be the modulus of a random polynomial near the unit circle.

In this work, the lcd allows us to develop clear arguments to estimate how close are the roots of a Kac polynomial to the unit circle. To be more precise, when the rcs of a Kac polynomial are iid rvs with zero mean and finite second moment, the majority of the roots are a distance of order O$( n^{-2}(\log{n})^{-3})$ with probability $1-\textnormal{O}((\log n)^{-{1}/{2}})$.  The main obstacle to extend this result comes from the Salem-Zygmund inequality as we will see in Section \ref{section2}.

The main result of this paper is the following.
\begin{theorem}\label{thm28052018}
Let $\left\{\xi_j: j\geq 0\right\}$ be a sequence of real iid non-degenerate real rvs satisfying
\begin{equation}
\label{H}
\tag{\bf H}
\sup_{u\in\R} \mathbb{P}(|\xi_0 - u| \leq \gamma) \leq 1-q\;\; \mbox{ and } \;\;\mathbb{P}(|\xi_0|>M) \leq \frac{q}{2} 
\end{equation}
for some  $M>0$, $\gamma>0$ and $q\in (0,1)$.
Assume that $\mathbb{E}[\xi_0] =0$  and $\mathbb{E}[\xi_0^2]<\infty$. Then for all $t\geq 1$ fixed we have 
\begin{equation}\label{eq:theorema}
\mathbb{P}\Big(\min_{\abs{\abs{z}-1} \leq tn^{-2}\lP\log n\rP^{-3}} \abs{W_n(z)} \leq t n^{-{1}/{2}} \lP\log n\rP^{-2}\Big) = \textnormal{O}((\log n)^{-{1}/{2}}),
\end{equation}
where the implicit constant depends on $t$ and the distribution of $\xi_0$.
\end{theorem}

\begin{remark}
~
\begin{enumerate}
\item In Theorem \ref{thm28052018} only is assumed  the finiteness of the second moment, zero mean 
and condition \eqref{H}
which include
Rademacher and standard Gaussian rvs.
As a directly consequence of Theorem \ref{thm28052018},  we have
\[
\mathbb{P}\big(W_n \textrm{ has no roots on }
\big\{z\in \mathbb{C}:\abs{\abs{z}-1} \leq tn^{-2}\lP\log n\rP^{-3}\big\}
\big)= 1-\textnormal{O}((\log n)^{-{1}/{2}}).
\]
\item We point out that in \eqref{eq:theorema}  we consider  the minimum of the modulus of the Kac polynomial over  the set
$\{\abs{\abs{z}-1} \leq tn^{-2}\lP\log n\rP^{-3}\}$ which is properly contained into the region considered in expression
\eqref{eq:kay}, but it contains the region considered in \eqref{eq:limitzero}.
Nevertheless, we 
obtain the upper bound
$\textnormal{O}((\log n)^{-{1}/{2}})$ which  improves the bound given in \eqref{eq:kay}. 
\end{enumerate}
\end{remark}

This manuscript is organized as follows.
In Section \ref{section2} we give an outline of the proof.
In Section \ref{section3} we provide the proof of Theorem \ref{thm28052018}. Finally, in Appendix \ref{A:app270120191754} we prove auxiliary results that we used throughout the manuscript.

\section{Outline of the proof}\label{section2}
In this section, we present the strategy used to prove Theorem \ref{thm28052018}.
Our goal is to estimate $\mathbb{P}(A_n)$, where
\begin{equation*} 
{A}_n:=\Big\{ \min_{z\in \mathbb{C}:\abs{\abs{z}-1} \leq tn^{-2}\lP\log n\rP^{-3}} \abs{W_n(z)} \leq  t n^{-{1}/{2}} (\log n)^{-2}\Big\}
\end{equation*}
and $t\geq 1$ is a fix constant. 
First, motivated by the estimates  given in \cite{Kon1999}, Section 2, p. 4964, we analyze the probability of the events
\begin{equation*}
A_{n,\alpha}:=\{ \abs{W_n\lP \exp( i2\pi x_\alpha)\rP} \leq g_n \} \quad \textrm{ for }
x_\alpha = \frac{\alpha}{N_n}, ~\alpha=0,\ldots, N_n-1,
\end{equation*} 
where $N_n$ and $g_n$ are appropriate functions of $n$ (later on, we provide precise description of them). 
We anticipate that $N_n\approx n^2(\log(n))^3$, which is similar to the number of balls using in \cite{Kon1999}.
We point out that $N_n$ needs to trade off $g_n$ in order to the probability of $A_{n,\alpha}$ tends to zero, as $n\to \infty$.

Second, for each $\alpha=0,\ldots,N_n-1$
we analyze the arithmetic structure of the sequence 
$\{\exp(i 2\pi j x_\alpha) : j=0,\ldots, n-1\}$
and using the so-called small ball inequalities we prove that $\mathbb{P}(A_{n,\alpha})\to 0$, as $n\to \infty$.
The idea is to apply the Taylor Theorem to 
approximate $W_n$ in  small balls with centers at $\exp( i2\pi x_\alpha)$. It allows us to
write the event $A_n$ as the union of events of the form $A_{n,\alpha}$.
However,  we need to handle the maximum value for the derivative of $W_n$ on the unit circle.
The latter can be done by a Salem-Zygmund type inequality, which consists in estimate the maximum possible value of a Kac polynomial on the unit circle.

Denote by $\|W^{\prime}_n\|_\infty$ the supremum norm of $W^{\prime}_n$ over the unit circle.  In the case of $\xi_0,\ldots,\xi_{n-1}$ being iid sub-Gaussian rvs,  a Salem-Zygmund type inequality (in probability) gives
\begin{equation}\label{eqn210220191805}
\Prob{\|W_n\|_\infty> C_p n^{{1}/{2}}\lP \log n\rP^{{1}/{2}}} = \textnormal{O}\lP {n^{-2}}\rP
\end{equation} 
for some suitable positive constant $C_p$. In \cite{GP}, the authors showed that \eqref{eqn210220191805} holds for iid 
zero mean
rvs with finite moment generating function. In this paper, {\it{we are not assuming the existence of the moment generating function}}. Instead of, we assume 
finiteness of the second moment.
By applying the majorizing measure method,
Weber \cite{Weber2006} showed 
\eqref{eqn210220191805} in expectation. 
To be more precise, let $\xi_0,\xi_1,\ldots,\xi_{n-1}$ be iid zero mean rvs with finite second moment, Corollary 2 in \cite{Weber2006} implies that there exists a positive constant 
$\tilde{C}$ (only depends on $\mathbb{E}[\xi_0^2]$) such that 
\[
\E{\|W_n\|_\infty} \leq \tilde{C} n^{{1}/{2}} \lP \log n\rP^{
{1}/{2}}\quad \textrm{  for any } n\in \mathbb{N}.
\]
To improve  Theorem \ref{thm28052018} (using lcd technique) to more general rcs
we require a refined version of the Salem-Zygmund inequality for rvs without finite second moment. At the moment, the authors are not able to obtain a Salem-Zygmund type inequality for rvs without the finiteness of the  second moment.
Later, we apply small ball inequalities to show that
\[
\Prob{ \abs{W_n\lP \exp\lP i2\pi x_\alpha\rP\rP} \leq g_n} \rightarrow 0,\quad\textrm{ as } n\to\infty.
\] 
This kind of inequalities allow us to consider more general rcs and provide a new proof of the main theorem in \cite{Kon1999}. 
To apply small ball inequalities we analyze the lcd for some specific matrix.
In the sequel, we give the definition of the lcd for a matrix. Denote
$\log_+ x:= \max\lL\log x, 0\rL$ for any $x>0$.

\begin{definition}[Least common denominator (lcd)]\label{def260920191725} Let $L>0$ be a positive number. Let $\|\cdot\|_2$ be the standard Euclidean norm and let $ \textnormal{dist}\lP v, \mathbb{Z}^M\rP$ denote the distance between the vector $v\in\R^M$ and the set $\mathbb{Z}^M$.   For a given matrix $V\in\R^{m\times M}$ the lcd is defined as 
\[
D(V) := \inf\lL \|{\Theta}\|_2 : \Theta\in\R^m, \textnormal{dist}\lP V^T\Theta, \Z^M\rP < L\sqrt{ \log_+\frac{\|{V^T\Theta}\|_2}{L}} \rL.
\]
\end{definition}
For a review of the concept lcd, we recommend Section 7 in \cite{RV2}. 
For our purposes, in Definition \ref{def260920191725} we take $m=2$, $M=n$ 
and the matrix $V$ is given by
\[ 
V:=\lC
\begin{array}{cccc}
1 & \cos\lP 2\pi x_\alpha \rP & \ldots & \cos\lP (n-1)2\pi x_\alpha \rP \\
0 & \sin\lP 2\pi x_\alpha \rP & \ldots & \sin\lP (n-1)2\pi x_\alpha \rP
\end{array}
\rC.
\] 
Set $X=\lC\xi_0,\ldots,\xi_{n-1}\rC^T$. Observe that
\[
\Prob{ \|{V X}\|_2 \leq g_n} = \Prob{ \abs{W_n\lP \exp\lP i2\pi x_\alpha\rP\rP} \leq g_n}.
\] 
Note that if $\det\lP VV^T\rP>0$, Theorem 7.5 in \cite{RV2} implies that for $a>0$ and $t\geq 0$ 
\begin{align}\label{rvi}
\Prob{ \|aV X\|_2 \leq t } \leq 
\frac{C^2 L^2}{2a^2(\det(VV^T))^{{1}/{2}}}
\lP t + \frac{1}{D(a V)}\rP^2,
\end{align}where $L\geq \sqrt{{2}/{q}}$ and the constant $C$ only depends on constants $M$, $\gamma$, $q$ specified in Theorem \ref{thm28052018}.
By Definition \ref{def260920191725} it is not hard to deduce that for any $a>0$, $D(aV)\geq  (\nicefrac{1}{a}) D(V)$. 
Recall the inequality $(x+y)^2\leq 2x^2 + 2y^2$ for any $x,y\in\R$.
By \eqref{rvi} we deduce that

\begin{align}
\label{eqn22022019901}
\Prob{ \|aV X\|_2 \leq t } & \leq  \frac{C^2L^2t^2}{a^2\lP\det\lP VV^T\rP\rP^{1/2}} + \frac{C^2L^2}{a^2\lP\det\lP VV^T\rP\rP^{1/2}\lP D(aV)\rP^2} \nonumber\\
 & \leq  \frac{C^2L^2t^2}{a^2\lP\det\lP VV^T\rP\rP^{1/2}} + \frac{C^2L^2}{\lP\det\lP VV^T\rP\rP^{1/2}\lP D(V)\rP^2}.\tag{\mbox{\bf K}}
\end{align}
Since $x_\alpha = \frac{\alpha}{N_n}$, the arithmetic properties of $x_\alpha$ given by $\alpha$ and $N_n$ should play an important role in the estimates. Depending on the greatest common divisor between $\alpha$ and $N_n$, $\gcd\lP \alpha, N_n\rP$, we deduce suitable positive lower bounds for $\det\lP V^T V\rP$ and $\textnormal{dist}\lP V^T \Theta, \Z^n\rP$ which together with \eqref{eqn22022019901} allow us to show that
$\Prob{ \|{V X}\|_2 \leq g_n}$ is sufficiently small.

\subsection*{Taylor's approximation}
In the sequel, 
define the trigonometric random polynomial $T_n(x):= \sum_{j=0}^{n-1} \xi_j e^{ijx}$,  $x\in\R$ and denote by $T^{\prime}_n$ its derivative with respect to $x$. To do the notation shorter, we
denote by $\Delta_n$ the following event:
\[
\Delta_n := \lL \max_{z\in \mathbb{C}:\abs{\abs{z}-1} \leq 2t n^{-11/10}} |W_n(z)| \leq n^{3/2} , \;\|T^{\prime}_n\|_\infty \leq C_0 n^{3/2} \log n\rL,
\]
where $C_0$ is a positive constant that we will precise later.
We also denote by $\Prob{A,B}$ the probability $\Prob{A\cap B}$ for any two events $A$ and $B$.
By the total probability law, we deduce

\begin{eqnarray}\label{eqn071020191038}	
& &\Prob{A_n}  \leq \;  \Prob{A_n, \Delta_n} 
	+ \Prob{\|{T^{\prime}_n}\|_\infty > C_0 n^{\frac{3}{2}} \log n} 
	\nonumber \\& &\hspace{2cm} 
	+\;\Prob{\max_{z\in \mathbb{C}:|\abs{z}-1| \leq  2t n^{-11/10}} |W_n(z)|  > n^{3/2}}.
\end{eqnarray}	
The Markov inequality yields 
\begin{eqnarray*}
& & \hspace{-0.5cm} \Prob{\max_{z\in \mathbb{C}:|\abs{z}-1| \leq 2t n^{-11/10}} |W_n(z)|  > n^{3/2}}  \leq  \Prob{\sum_{j=0}^{n-1} \abs{\xi_j} \lP 1+ \frac{2t}{n^{1+1/10}}\rP^j> n^{3/2}} \nonumber \\
& &\hspace{2.5cm} \leq  \frac{\E{\sum_{j=0}^{n-1} \abs{\xi_j} \lP 1+ \frac{2t}{n^{1+1/10}}\rP^j }}{n^{3/2}}
\leq  \frac{e^{2t} n \E{|\xi_0|}}{n^{3/2}} = \frac{e^{2t} \E{|\xi_0|}}{n^{1/2}}.\nonumber\\
\end{eqnarray*}
In other words, we have
\begin{equation}\label{eqn07102019}
\Prob{ \max_{z\in \mathbb{C}:\abs{\abs{z}-1} \leq 2t n^{-11/10}} \abs{W_n(z)}  > n^{3/2}}=\textnormal{O}\lP n^{-1/2}\rP,
\end{equation} where the implicit constant depends on $t$ and $\E{|\xi_0|}$.
On the other hand, the Bernstein inequality (see Theorem 14.1.1 in \cite{RahSch2002}) allows us to deduce 
 for the second term in the right-side hand of \eqref{eqn071020191038} that
\[
\Prob{\|T^{\prime}_n\|_\infty > C_0 n^{3/2} \log n} \leq  \Prob{\|{T_n}\|_\infty > C_0 n^{1/2} \log n}.
\]
Since $\E{\xi_0}=0$ and $\E{\xi^2_0}<\infty$, one can apply Corollary 2 of \cite{Weber2006} which together with the Markov inequality imply
\[
\begin{split}
\Prob{\|T_n\|_\infty > C_0 n^{1/2} \log n}  \leq 
\frac{C(\E{\xi^2_0})^{1/2} n^{1/2} \lP\log n\rP^{1/2}}{C_0 n^{1/2} \log n} = \frac{C(\E{\xi^2_0})^{1/2}}{C_0 \lP\log n\rP^{1/2}},
\end{split}
\]
where $C$ is a universal positive constant. 
Consequently,  the Bernstein inequality yields
\begin{equation}\label{eqn071020191021}
\Prob{\|T^{\prime}_n\|_\infty > C_0 n^{3/2} \log n}  = \textnormal{O}\lP{\lP\log n\rP^{-{1}/{2}}}\rP.
\end{equation} 
By inequalities \eqref{eqn071020191038}, \eqref{eqn07102019} and \eqref{eqn071020191021}, we observe 
that to estimate $\mathbb{P}(A_n)$ we only need to analyze $\mathbb{P}(A_n,\Delta_n)$.
\begin{remark}
In the preceding reasoning we only use zero mean and finiteness of the second moment of $\xi_0$.
In particular, 
it holds for sub-Gaussian rvs which includes Rademacher, standard Gaussian, and bounded rvs.
\end{remark}

\subsection*{Arithmetic properties of $\mathbf{x_\alpha}$}
In the sequel, we decompose the event $A_n\cap \Delta_n$ into regions for which the arithmetic properties of $x_\alpha$ 
are useful in order to use the anti-concentration assumption \eqref{H} and allows us to show 
 $\mathbb{P}(A_n,\Delta_n)$ tends to zero, as $n\to \infty$.
We point out that in the following reasoning we only use assumption \eqref{H}.

To achieve our goal, we consider a set of balls with center at points on the unit circle with an adequate radius. We distinguish two kind of balls. The special balls with center at $1+0i$ and $-1+0i$, where the radius $r$ is {\it large} 
($r=2tn^{-11/10}$), and the balls with center at points $z$ with argument satisfying $n^{-11/10}<\abs{\arg(z) } < \pi - n^{-11/10} $ and {\it small} radius ($r=2t n^{-2}\lP\log n\rP^{-3}$).

Recall that for any $x\in \mathbb{R}$, $\left\lfloor x \right\rfloor$  denotes the greatest integer less than or equal to $x$.
Let $N= \lfloor n^2\lP \log n \rP^3 \rfloor$ and $x_\alpha= \frac{\alpha}{N}$ for $\alpha=0,1,\ldots,N-1$. For $a\in\C$ and $s>0$, denote by $\textrm{B}\lP a, s\rP$ the closed ball with center $a$ and radius $s$, i.e., $\textrm{B}\lP a, s\rP = \lL z\in\C : \abs{z-a} \leq  s\rL$. Denote by 
$\mathbb{S}^1$ the unit circle.
Let 
\begin{align*}
\mathcal{A}\left(\mathbb{S}^1,tn^{-2}\lP\log n\rP^{-3}\right):=\left\{z\in \mathbb{C}:\abs{\abs{z}-1} \leq tn^{-2}\lP\log n\rP^{-3}\right\}.
\end{align*}
Notice that
\begin{align*} 
 \mathcal{A}\left(\mathbb{S}^1,tn^{-2}\lP\log n\rP^{-3}\right) = & 
\lL z\in\mathcal{A} : n^{-11/10}<\abs{\arg(z) } < \pi - n^{-11/10} \rL \\ 
&\hspace{-1cm}\cup \lL z\in\mathcal{A} : \abs{\arg(z)}\leq n^{-11/10}\quad \mbox{ or }\quad \abs{\arg(z)-\pi} \leq n^{-11/10} \rL.
\end{align*}
Let $t\geq 1$ and observe that
\begin{align*} 
&\lL z\in\mathcal{A} : \abs{\arg(z)}\leq n^{-11/10}\quad \mbox{ or }\quad \abs{\arg(z)-\pi} \leq n^{-11/10} \rL \\
&\hspace{3cm} \subset \B\lP -1+0i, 2tn^{-11/10}\rP \cup \B\lP 1+0i, 2tn^{-11/10}\rP.
\end{align*} 
The preceding inclusion yields that any $z\in\mathcal{A}$ with \textit{small argument} belongs in the union of the  balls with center at $1+0i$ and $-1+0i$ with radius $2tn^{-11/10}$. On the other hand, for $z\in\mathcal{A}$ with \textit{large argument} we have
\begin{align*}
&\lL z\in\mathcal{A} : n^{-11/10}<\abs{\arg(z) } < \pi - n^{-11/10} \rL \\
&\hspace{2cm}\subset \bigcup^{N-1} _{\substack{\alpha=1 \\ \alpha\;:\; n^{-11/10}<\abs{2\pi x_\alpha} < \pi - n^{-11/10}}}
\B \lP e^{i2\pi x_\alpha},2tn^{-2}\lP\log n\rP^{-3}\rP.
\end{align*}
Define
\[
\begin{split}
J_1(n,N)&:=\left\{\alpha\in [1,N-1]\cap \mathbb{N}: \gcd\lP \alpha, N\rP \geq n^{11/10} \lP\log n\rP^{-1/2}\right\},\\
J_2(n,N)&:=\left\{\alpha\in [1,N-1]\cap \mathbb{N}: n^{11/10} \lP\log n\rP^{-1/2}\geq \gcd\lP\alpha,N\rP\geq n\lP\log n\rP^3\right\},\\
J_3(n,N)&:=\left\{\alpha\in [1,N-1]\cap \mathbb{N}: n\lP\log n\rP^3 \geq \gcd\lP \alpha, N\rP \geq n^{9/10}\lP \log n\rP^3\right\},
\end{split}
\]
where $\gcd(\alpha,N)$ denotes the greatest common divisor between $\alpha$ and $N$.
For any $\alpha\in J_3(n,N)$ we have
\[
n - \frac{1}{n\lP \log n\rP^3}\leq  \frac{N}{\gcd\lP \alpha, N\rP} \leq n^{11/10}.
\] 
The preceding inequalities mean that the irreducible fraction of $x_\alpha$ is as small as a multiple of $n^{-11/10}$. Therefore,
\begin{align*}
&\bigcup^{N-1} _{\substack{\alpha=1 \\ \alpha\;:\; n^{-11/10}<\abs{2\pi x_\alpha} < \pi - n^{-11/10}}}
\B \lP e^{i2\pi x_\alpha},2tn^{-2}\lP\log n\rP^{-3}\rP\\
&\hspace{0.5cm}= \bigcup_{\alpha\in J_1(n,N)}
\B \lP e^{i2\pi x_\alpha},2tn^{-2}\lP\log n\rP^{-3}\rP
\cup 
\bigcup_{\alpha\in J_2(n,N)}
\B \lP e^{i2\pi x_\alpha},2tn^{-2}\lP\log n\rP^{-3}\rP\\
&\hspace{1cm} \cup 
\bigcup_{\alpha\in J_3(n,N)}
\B \lP e^{i2\pi x_\alpha},2tn^{-2}\lP\log n\rP^{-3}\rP
\end{align*}
We emphasize that if $\alpha\in J_1(n,N)\cup J_2(n,N)\cup J_3(n,N)$. Then we have 
\[n^{-11/10} < \abs{2\pi x_\alpha} < \pi - n^{-11/10}.\] 
Consequently, 
\begin{align}
\Prob{{A}_n, \Delta_n} \; \leq & \; \Prob{\Delta_n, \min_{z\in \B\lP 1+0i, 2tn^{-11/10}\rP} \abs{W_n(z)} < tn^{-1/2}\lP \log n\rP^{-2}} \nonumber\\
& +\Prob{\Delta_n, \min_{z\in \B\lP -1 + 0i, 2tn^{-11/10}\rP} \abs{W_n(z)} < tn^{-1/2}\lP \log n\rP^{-2}} \label{mainestimate} \\
&  + \sum\limits_{\alpha\in J_1(n,N)} \Prob{\Delta_n, \B_\alpha}+
\sum\limits_{\alpha\in J_2(n,N)} \Prob{\Delta_n, \B_\alpha}+\sum\limits_{\alpha\in J_3(n,N)} \Prob{\Delta_n, \B_\alpha},\nonumber
\end{align}
where
\[
\B_\alpha := \lL \min_{z\in \B\lP e^{i2\pi x_\alpha}, 2t n^{-2}\lP \log n\rP^{-3}\rP} \abs{W_n(z)} <tn^{-1/2}\lP \log n \rP^{-2} \rL.
\]
The right-hand side of \eqref{mainestimate} will be estimated as follows.
\begin{lemma}\label{lem:A}
The following holds.
\begin{equation*}
\Prob{\Delta_n, \min_{z\in \B\lP 1+0i, 2tn^{-{11}/{10}}\rP} \abs{W_n(z)} < tn^{-{1}/{2}}\lP \log n\rP^{-2}} =\textnormal{O}\lP\frac{\log n}{n^{{1}/{10}}}\rP
\end{equation*} 
and 
\begin{equation*}
\Prob{\Delta_n, \min_{z\in \B\lP -1+0i, 2tn^{-11/10}\rP} \abs{W_n(z)} < tn^{-1/2}\lP \log n\rP^{-2}} =\textnormal{O}\lP\frac{\log n}{n^{{1}/{10}}}\rP,
\end{equation*}
where  
the implicit constants in the 
big $\textnormal{O}$-notation 
depend on $L$ and $t$. 
\end{lemma}

\begin{lemma}\label{lem:B}
Assume that $\gcd\lP \alpha, N\rP \geq n^{11/10} \lP\log n\rP^{-1/2}$, where
$N= \lfloor n^2( \log n )^3 \rfloor$. Then
 for a suitable constant $\tilde{C}$ it follows
\begin{equation*}\label{eqn071020191626}
\sum\limits_{\alpha \in J_1(n,N)}  \Prob{\abs{W_n\lP e^{i2\pi x_\alpha} \rP} \leq \tilde{C} t n^{-1/2}\lP\log n\rP^{-2}} = \textnormal{O}\lP \frac{\lP\log n\rP^{4}}{n^{1/20}}\rP,
\end{equation*} 
where
the implicit constant 
in the 
big $\textnormal{O}$-notation  
 depends on $L$ and $t$. 
\end{lemma}

\begin{lemma}\label{lem:C}
Assume that $\frac{n^{11/10}}{\lP \log n\rP^{1/2}}\geq \gcd\lP\alpha,N\rP\geq n\lP\log n\rP^3$, 
where \\
$N= \lfloor n^2( \log n )^3 \rfloor$. Then
 for a suitable constant $\tilde{C}$ it follows
\begin{equation*}\label{eqn071020191807}
\sum\limits_{\alpha \in J_2(n,N)}  \Prob{\abs{W_n\lP e^{i2\pi x_\alpha} \rP} \leq \tilde{C}t  n^{-1/2}\lP\log n\rP^{-2}} = \textnormal{O}\lP\frac{1}{\log n}\rP,
\end{equation*} 
where
 the implicit constant
in the 
big $\textnormal{O}$-notation  
  depends on $L$ and $t$. 
\end{lemma}

\begin{lemma}\label{lem:D}
Assume that $n\lP\log n\rP^3 \geq \gcd\lP \alpha, N\rP \geq n^{9/10}\lP \log n\rP^3$, where \\
$N= \lfloor n^2\lP \log n \rP^3 \rfloor$. Then
 for a suitable constant $\tilde{C}$ it follows
\begin{equation*}\label{eqn071020191749}
\sum\limits_{\alpha \in J_3(n,N)}  \Prob{\abs{W_n\lP e^{i2\pi x_\alpha} \rP} \leq \tilde{C} t  n^{-1/2}\lP\log n\rP^{-2}} = \textnormal{O}\lP\frac{1}{n^{1/10}}\rP, 
\end{equation*} 
where
 the implicit constant 
 in the 
big $\textnormal{O}$-notation 
depends on $L$ and $t$. 
\end{lemma}

In the sequel, we stress the fact that Theorem \ref{thm28052018}
is just a consequence
of what we have already stated up to here.
Indeed, 
combining  
Lemma \ref{lem:A}, Lemma \ref{lem:B}, Lemma \ref{lem:C}, Lemma \ref{lem:D}, estimate
\eqref{eqn07102019} and estimate \eqref{eqn071020191021}
in inequality \eqref{eqn071020191038}
yield Theorem \ref{thm28052018}.
 
\section{Proof of Theorem \ref{thm28052018}}\label{section3}
In this section, we show that the left-hand side of inequality \eqref{mainestimate} is of order 
$\textrm{O}((\log(n))^{-{1}/{2}})$.
\subsection{Estimates on the balls centered at $-1$ and $1$}
\begin{proof}[Proof of Lemma \ref{lem:A}]
Let $z\in \B\lP 1+0i, 2tn^{-11/10}\rP$. 
The Taylor Theorem implies
\[
\abs{W_n(z) - W_n(1)} \leq \abs{z - 1} \abs{W^{\prime}_n(1)} + \abs{R_2(z)},
\] 
where $R_2(z)$ is the error of the Taylor approximation of order $2$. On $\Delta_n$  we have
\[
\begin{split}
\abs{R_2(z)} &  \leq    \; \frac{\lP 2 tn^{-1-1/10}\rP^2}{1-\textnormal{o}(1)}  \max_{z\in \B\lP 1+0i, 2tn^{-11/10}\rP} \abs{W_n(z)} \\
& \leq \; \frac{4t^2 n^{-2-1/5} n^{3/2}}{1-\textnormal{o}(1)} = \frac{4t^2 n^{-1/2-1/5}}{1-\textnormal{o}(1)},
\end{split}
\] where $\textnormal{o}(1) = 2 tn^{-1-1/10}$. The preceding inequality and assuming that $\Delta_n$ holds yield
\begin{eqnarray*}
\abs{W_n(z) - W_n(1)} & \leq & 2tn^{-1-1/10} \abs{W^{\prime}_n(1)} + \frac{4t^2 n^{-1/2-1/5}}{1-\textnormal{o}(1)} \\
& \leq & 2tn^{-1-1/10} \|T^{\prime}_n\|_\infty + \frac{4t^2 n^{-1/2-1/5}}{1-\textnormal{o}(1)} \\
& \leq & 2C_0 t n^{1/2-1/10} \log n + \frac{4t^2 n^{-1/2-1/5}}{1-\textnormal{o}(1)}.
\end{eqnarray*} 
Hence,
\[
\begin{split}
\Prob {\Delta_n, \min_{z\in \B\lP 1+0i, 2tn^{-11/10}\rP} \abs{W_n(z)} < tn^{-1/2}\lP \log n\rP^{-2}} & \leq  \\
& \hspace{-3cm}\Prob{\abs{W_n(1)} \leq 2C_2 t n^{1/2-1/10}\log n},
\end{split}
\]
where $2C_2 = 2C_0 + 4t + 1$. 
Since $W_n(1)=\sum_{j=0}^{n-1} \xi_j$, Corollary 7.6 in \cite{RV2} implies
\[
\Prob{\abs{W_n(1)} \leq 2C_2 t n^{1/2-1/10}\log n} \leq \frac{C_3 L}{\|{\mathbf{a}}\|}\lP 2C_2t + \frac{1}{D(\mathbf{a})}\rP
\]
for $L\geq \sqrt{1/q}$,
where $C_3$ is a positive constant and  $D(\mathbf{a})$ is the lcd of 
\[
\mathbf{a}=\left(n^{1/2-1/10} \log n\right)^{-1} (1,\ldots,1) \in \R^n.
\]
By Proposition 7.4 in \cite{RV2} we have $D(\mathbf{a})\geq {1}/{2} n^{1/2 - 1/10}\log n$. Then
\[
\begin{split}
\Prob{\abs{W_n(1)} \leq 2C_2 t n^{1/2-1/10}\log n} & \leq \; \frac{C_3 L \log n}{n^{1/10}}\lP 2C_2 t + \frac{2}{n^{1/2-1/10} \log n}\rP \\
& \leq \; \frac{\lP2C_2t+2\rP C_3 L \log n}{n^{1/10}}.
\end{split}
\]
Therefore,
\begin{equation*}\label{eqn071020191625}
\Prob{\Delta_n, \min_{z\in \B\lP 1+0i, 2tn^{-11/10}\rP} \abs{W_n(z)} < tn^{-1/2}\lP \log n\rP^{-2}} =\textnormal{O}\lP\frac{\log n}{n^{1/10}}\rP.
\end{equation*} 
On the other hand, for $z\in B\lP -1+0i, 2tn^{-11/10}\rP$ similar reasoning yields
\[
\begin{split}
&\Prob{\Delta_n, \min_{z\in \B\lP -1+0i, 2tn^{-11/10}\rP} \abs{W_n(z)} < tn^{-1/2}\lP \log n\rP^{-2}}  \leq  \\
& \hspace{6cm} \Prob{\abs{W_n(-1)} \leq 2C_2 t n^{1/2-1/10}\log n}.
\end{split}
\]
In this case, we need to analyze $W_n(-1)=\sum_{j=0}^{n-1} \lP -1\rP^{j} \xi_j$.
Again taking $L\geq \sqrt{1/q}$ and applying Corollary 7.6 in \cite{RV2} we obtain
\begin{eqnarray*}
\Prob{\abs{W_n(-1)} \leq 2C_2 t n^{1/2-1/10}\log n} \leq \frac{C_3 L}{\|{\mathbf{b}}\|}\lP 2C_2t + \frac{1}{D(\mathbf{b})}\rP,
\end{eqnarray*} 
where $C_3$ is a positive constant and $D(\mathbf{b})$ is the lcd of 
\[\mathbf{b}=\lP n^{1/2-1/10}\log n\rP^{-1}\lP1,-1,1,\ldots,(-1)^{n-1}\rP\in\R^n.\] 
By Proposition 7.4 in \cite{RV2} we have $D(\mathbf{b})\geq {1}/{2} n^{1/2-1/10}\log n$. Then
\[
\begin{split}
\Prob{\abs{W_n(-1)} \leq 2C_2 t n^{1/2-1/10}\log n} \; & \leq \;  \frac{C_3 L \log n}{n^{1/10}}\lP 2C_2 t + \frac{2}{n^{1/2-1/10} \log n}\rP \\
   &  \leq  \; \frac{\lP2C_2t+2\rP C_3 L \log n}{n^{1/10}}.
\end{split}
\]
Therefore, 
\begin{equation*}\label{eqn071020191351}
\Prob{\Delta_n, \min_{z\in \B\lP -1+0i, 2tn^{-11/10}\rP} \abs{W_n(z)} < tn^{-1/2}\lP \log n\rP^{-2}}   = \textnormal{O}\lP\frac{\log n}{n^{1/10}}\rP.
\end{equation*} 
\end{proof}

\subsection{Estimates of $\Prob{\Delta_n, \B_\alpha }$}
In the sequel, we apply the Taylor Theorem repeatedly in order to reduce $\Prob{ \Delta_n, \B_\alpha }$ to an estimate of the probability of how small a sum of iid rvs can be. The latter can be computed (estimated) using small ball inequalities.

Let  $z\in \B\lP e^{i2\pi x_\alpha}, 2t n^{-2} \lP\log n \rP^{-3}\rP$ and  assume that $\Delta_n$ holds. The Taylor Theorem yields
\[
\begin{split}
\abs{ W_n(z) - W_n\lP e^{i2\pi x_\alpha} \rP} \; \leq  \; & \abs{z - e^{i2\pi x_\alpha}} \abs{W^{\prime}_n \lP e^{i2\pi x_\alpha} \rP} + \abs{R_2(z)} \\
\; \leq \; &  2tn^{-2}\lP \log n\rP^{-3} \abs{W^{\prime}_n \lP e^{i2\pi x_\alpha} \rP} + \frac{4t^2 n^{-5/2}\lP\log n\rP^{-6}}{1-\textnormal{o}(1)} \\
\; \leq  \;& \lP  2tC_0 +4t^2\rP n^{-1/2} \lP \log n\rP^{-2},
\end{split}
\] 
where $\textnormal{o}(1)=2tn^{-2} \lP\log n\rP^{-3}$.
Then
\begin{equation}\label{eqn071020191353}
\Prob{\Delta_n,\B_\alpha} \leq \Prob{\abs{W_n\lP e^{i2\pi x_\alpha} \rP} \leq 2tC_2 n^{-1/2}\lP\log n\rP^{-2}}.
\end{equation}
To show that $\Prob{\Delta_n,\B_\alpha}$ tends to zero as $n$ increases, we rewrite the sum $W_n\lP e^{i2\pi x_\alpha}\rP$ as the product of a matrix by a vector, and then we analyze the lcd of the corresponding matrix.

Define the $2\times n$ matrix $V_\alpha$ as follows
\[
V_\alpha := 
\lC
\begin{array}{cccc}
1 & \cos\lP 2\pi x_\alpha \rP & \ldots & \cos\lP (n-1)2\pi x_\alpha \rP \\
0 & \sin\lP 2\pi x_\alpha \rP & \ldots & \sin\lP (n-1)2\pi x_\alpha \rP
\end{array}
\rC
\] 
and take $X=\lC\xi_0,\ldots,\xi_{n-1}\rC^T\in\R^n$. 
Notice that
\[
\|{V_\alpha X }\|_2 = \Big|\sum_{j=0}^{n-1} \xi_j e^{ij2\pi x_\alpha}\Big| = \abs{W_n\lP e^{i2\pi x_\alpha}\rP}.
\]
Let $\Theta=r \lC \cos\lP\theta\rP, \sin\lP\theta\rP \rC^T\in\R^2$, where  $r>0$ and $\theta\in[0,2\pi)$. For fix $r$ and $\theta$ we have
\[
V_\alpha^T \Theta = r \lC \cos\lP  - \theta\rP, \cos\lP 2\pi x_\alpha - \theta\rP, \ldots, \cos\lP 2(n-1)\pi x_\alpha - \theta\rP \rC^T\in \R^n.
\]
We also point out that   $\|V_\alpha^T \Theta \|_2\leq r\sqrt{n}$.
On the other hand, we observe that
\[
\det\lP V_\alpha V_\alpha^T \rP = \det\lC
\def\arraystretch{1.5}
\begin{array}{cc}
\sum_{j=0}^{n-1} \cos^2\lP j2\pi x_\alpha \rP & \frac{1}{2}\sum_{j=0}^{n-1} \sin\lP 2\cdot j2\pi x_\alpha \rP \\
\frac{1}{2}\sum_{j=0}^{n-1} \sin\lP 2\cdot j2\pi x_\alpha \rP & \sum_{j=0}^{n-1} \sin^2\lP j2\pi x_\alpha \rP
\end{array}
\rC.
\]
Bearing all this in mind, we can use the notion of lcd for high dimensions to obtain an accurate upper bound of the left-hand side of \eqref{eqn071020191353}. 

We recall that the events $\Delta_n\cap \B_\alpha$ are defined for  
\[n^{-11/10} < \abs{2\pi x_\alpha} < \pi - n^{-11/10}.\] 
Therefore, 
to estimate the left-hand side of \eqref{eqn071020191353} we distinguish the following three cases. 

\subsubsection{\textbf{Estimates on $J_1(n,N)$}} 
\begin{proof}[Proof of Lemma \ref{lem:B}]
Notice that
\[
\frac{N}{\gcd\lP \alpha, N\rP} \leq \frac{n^2\lP \log n\rP^{3}}{n^{11/10} \lP\log n\rP^{-1/2}} = n^{9/10} \lP \log n\rP^{7/2}.
\] 
and
\[
|2\pi x_{\alpha}|=2\pi \frac{\alpha}{N}=2\pi \frac{{\alpha}/{\gcd\lP \alpha, N\rP}}{{N}/{\gcd\lP \alpha, N\rP}}\geq 
2\pi \frac{1}{n^{9/10} \lP \log n\rP^{7/2}}.
\]
Then $2\pi x_\alpha$ also satisfies $n^{-1} < \abs{2\pi x_\alpha} < \pi - n^{-1}$ for all large $n$.
By Lemma 3.2 Part 1  in \cite{Kon1999}, there exist positive constants $c_4,C_4$  such that  
\begin{equation}\label{eqn110320191604}
c_4 n^2 \leq \det\lP V_\alpha V_\alpha^T\rP\leq C_4 n^2.
\end{equation}
By Lemma \ref{app270120191754} in Appendix \ref{A:app270120191754} we obtain that the number of indices $\alpha\in [1,N]\cap \N$ that satisfies the condition $\gcd\lP \alpha, N\rP \geq n^{11/10} \lP\log n\rP^{-1/2}$ is at most $\frac{N^{1+\textnormal{o}\lP 1 \rP}}{ n^{11/10}\lP\log n\rP^{-1/2}}$. By the definition of $N$ we obtain
\begin{equation}\label{eqn170920191648}
\frac{N^{1+\textnormal{o}\lP 1 \rP}}{ n^{11/10}\lP\log n\rP^{-1/2}} \leq \frac{n^{2+\textnormal{o}\lP 1 \rP} \lP \log n\rP^{7/2+\textnormal{o}\lP 1 \rP}}{n^{11/10}} = n^{9/10+\textnormal{o}\lP 1 \rP} \lP \log n\rP^{7/2+\textnormal{o}\lP 1 \rP}.
\end{equation}
By Proposition 7.4 in \cite{RV2}, the lcd of $V_\alpha$ satisfies $D\lP V_\alpha \rP \geq {1}/{2}$. Therefore, inequality \eqref{eqn22022019901}, inequality \eqref{eqn110320191604} and inequality \eqref{eqn170920191648} yield
\[
\begin{split}
 \sum\limits_{\alpha \in J_1(n,N)} & \Prob{\abs{W_n\lP e^{i2\pi x_\alpha} \rP} \leq 2tC_2 n^{-1/2}\lP\log n\rP^{-2}} \\ 
& \hspace{-1cm}
 \leq   n^{9/10+\textnormal{o}\lP 1 \rP} \lP \log n\rP^{7/2+\textnormal{o}\lP 1 \rP}\lP \frac{2C^2L^2\lP 2tC_2\rP^2}{\lP c_4n^2 \rP^{1/2} \lP n^{1/2} \lP\log n\rP^2\rP^2} + \frac{2C^2L^2}{\frac{1}{4} \lP c_4n^2 \rP^{1/2}}\rP \\ 
&\hspace{-1cm} \leq  \frac{8C^2C^2_2L^2t^2}{c_4^{1/2} n^{11/10-\textnormal{o}\lP 1 \rP}\lP\log n\rP^{1/2- \textnormal{o}\lP 1  \rP}} + \frac{8C^2L^2 \lP \log n\rP^{7/2+\textnormal{o}\lP 1 \rP}}{c_4^{1/2} n^{1/10-\textnormal{o}\lP 1 \rP}}
\end{split}
\]
for all large $n$.
Consequently,
\begin{equation*} 
\sum\limits_{\alpha \in J_1(n,N)} \Prob{\abs{W_n\lP e^{i2\pi x_\alpha} \rP} \leq 2tC_2 n^{-1/2}\lP\log n\rP^{-2}} = \textnormal{O}\lP \frac{\lP\log n\rP^{4}}{n^{1/20}}\rP,
\end{equation*} where the implicit constant depends on $L$ and $t$. 
\end{proof}

\subsubsection{\textbf{Estimates on $J_2(n,N)$}} 
\begin{proof}[Proof of Lemma \ref{lem:C}]
Notice that
\begin{equation}\label{arriba} 
n\geq \frac{N}{\gcd\lP \alpha,N\rP}\geq n^{9/10}\lP\log n\rP^{7/2} - 
\textnormal{o}(1),
\end{equation}
where $\textnormal{o}(1)=\frac{\lP \log n \rP^{1/2}}{n^{11/10}}$. The latter implies that
\[
|2\pi x_{\alpha}|=2\pi \frac{\alpha}{N}=2\pi \frac{{\alpha}/{\gcd\lP \alpha, N\rP}}{{N}/{\gcd\lP \alpha, N\rP}}\geq 
2\pi \frac{1}{n}.
\]
Then $2\pi x_\alpha$ also satisfies $n^{-1}\leq \abs{2\pi x_\alpha} \leq \pi - n^{-1}$ for all large $n$. By Lemma 3.2 Part 1 in \cite{Kon1999} there exist positive constants $c_4,C_4$ such that
\begin{equation}\label{eqn110320191606}
c_4 n^2 \leq \det\lP V_\alpha V_\alpha^T\rP \leq C_4 n^2.
\end{equation}
Note $x_\alpha=\frac{\alpha}{N}=\frac{\alpha'}{N'}$, where $\alpha= \alpha' \gcd\lP \alpha , N \rP$ and $N=N' \gcd\lP\alpha, N\rP$. Observe that  $\gcd\lP\alpha' , N'\rP = 1$. Since $N'\leq n$, for any $\theta$ we have
\[
\begin{split}
\lL \exp\lP i \lP j2\pi \frac{\alpha'}{N'} - \theta \rP\rP : j = 0,\ldots, N'-1 \rL = &\\
& \hspace{-3.5cm} \lL \exp\lP i \lP j2\pi \frac{1}{N'} - \theta \rP\rP : j = 0,\ldots, N'-1 \rL.
\end{split}
\]
Hence, without loss of generality, we assume that $x_\alpha = \frac{1}{N'}$.
A straightforward computation yields
\[
V_\alpha^T \Theta = r \lC \cos\lP  - \theta\rP, \cos\lP 2\pi x_\alpha - \theta\rP, \ldots, \cos\lP 2(n-1)\pi x_\alpha - \theta\rP \rC^T\in \R^n.
\] 
Notice that in the proof of Lemma \ref{lowerbound} in Appendix \ref{A:app270120191754} holds true for any real positive number $r$.
If $r\leq \frac{1}{32\pi x_\alpha}$, by Lemma~\ref{lowerbound} in Appendix \ref{A:app270120191754}, inequality \eqref{arriba} and remembering that $\|V_\alpha^T \Theta\|_2\leq r\sqrt{n}$  we deduce
\[
\begin{split}
\frac{1}{128\pi}\left(n^{9/10} \lP \log n\rP^{7/2} - \textnormal{o}\lP 1\rP\right) & \leq  \frac{1}{128\pi x_\alpha} 
\leq  \textnormal{dist}\lP V_\alpha^T\Theta, \Z^n\rP \\
&\hspace{-4cm} \leq  L\sqrt{\log_+\frac{\|{V_\alpha^T \Theta}\|_2}{L}} 
 \leq  L\sqrt{\log_+\frac{rn^{1/2}}{L}} \leq L\sqrt{\log_+\frac{n^{3/2} }{L}},
\end{split}
\]
which yields a contradiction due to $L\geq \sqrt{2/q}$ is fixed. Then  for $r>\frac{1}{32\pi x_\alpha}$ we have
\[
D\lP V_\alpha\rP\geq r>\frac{1}{32\pi}\left(n^{9/10} \lP \log n\rP^{7/2} - \textnormal{o}\lP 1\rP\right).
\]
Therefore, the preceding inequality together with inequality \eqref{eqn22022019901}, inequality   
\eqref{eqn110320191606} and the fact the cardinality of 
$J_2(n,N)$ is at most $N$ 
allow us to deduce
\[
\begin{split}
\sum\limits_{\alpha \in J_2(n,N)} & \Prob{\abs{W_n\lP e^{i2\pi x_\alpha} \rP} \leq 2tC_2 n^{-1/2}\lP\log n\rP^{-2}} \\ 
 & \hspace{-2cm} \leq n^2 \lP \log n\rP^3 \lP \frac{2C^2L^2 \lP2tC_2\rP^2}{\lP c_4 n^2\rP^{1/2} \lP n^{1/2} \lP\log n \rP^2\rP^2} \rP \\ 
 & \hspace{-1.8cm} \;\;+ n^2 \lP \log n\rP^3 \lP \frac{2C^2L^2}{\lP c_4 n^2\rP^{1/2} 
 \lP \frac{1}{32\pi}\left(n^{9/10} \lP \log n\rP^{7/2} - \textnormal{o}\lP 1\rP\right)\rP^2} \rP  \\ 
 & \hspace{-2cm}\leq \frac{8C^2C^2_2L^2t^2}{c^{1/2}_4\log n} + \frac{2048c\pi^2C^2 L^2}{c_4^{1/2}n^{2/5}\lP\log n\rP^4}
\end{split}
\]
for all large $n$,
where $c_4$ is a positive constant.
As a consequence we obtain
\begin{equation*}
\sum\limits_{\alpha \in J_2(n,N)}  \Prob{\abs{W_n\lP e^{i2\pi x_\alpha} \rP} \leq 2tC_2 n^{-1/2}\lP\log n\rP^{-2}} = \textnormal{O}\lP\frac{1}{\log n}\rP,
\end{equation*} where the implicit constant depends on $L$ and $t$.
\end{proof}

\subsubsection{\textbf{Estimates on $J_3(n,N)$}}
\begin{proof}[Proof of Lemma \ref{lem:D}]
This case requires a more refined analysis. Observe that
\[
n^{11/10} \geq \frac{N}{\gcd\lP \alpha, N\rP} \geq n - \textnormal{o}\lP 1\rP,
\]
where $\textnormal{o}\lP 1\rP=\frac{1}{n^2(\log n)^3}$.
Then  $2\pi x_\alpha$ satisfies
\[
n^{-11/10} \leq \abs{2\pi x_\alpha} \leq \lP n - \textnormal{o}\lP1\rP\rP^{-1}\quad \textrm{ or } \quad
\pi - \lP n - \textnormal{o}\lP1\rP\rP^{-1} \leq \abs{2\pi x_\alpha} \leq \pi - n^{-11/10}.
\] 
By Lemma 3.2 Part 2 in \cite{Kon1999}, there exist positive constants $c_4,C_4$ such that 
\begin{equation}\label{eqn110320191612}
c_4 n^{2-1/5} \leq \det\lP V_\alpha V_\alpha^T \rP \leq C_4 n^2.
\end{equation}
By Lemma \ref{app270120191754} we have that 
the number of indexes $\alpha\in [1,N]\cap \N$ that satisfy the condition $n\lP\log n\rP^3 \geq \gcd\lP \alpha, N\rP \geq n^{9/10}\lP \log n\rP^3$ is at most
$n^{11/10+\textnormal{o}(1)}(\log(n))^{\textnormal{o}(1)}$, where
$\textnormal{o}(1)\to 0$, as $n\to \infty$.

In the sequel, we analyze the lcd of $V_\alpha$. In particular, we find an appropriate lower bound for the distance between $V^T_\alpha\Theta$ and the set $\Z^n$. 
Since $x_\alpha=\frac{\alpha}{N}=\frac{\alpha'}{N'}$ with $\gcd\lP \alpha',N'\rP = 1$ and $N'\geq n - 1$ for all large $n$, then we have that all the points in $\lL \exp\lP i \lP j 2\pi x_\alpha - \theta \rP\rP : j=0,\ldots, n-1\rL$ are different between them.
Let $r\in \N$ consider the set of intervals of the form $\lC \frac{m}{r}, \frac{m+1}{r}\rC$ for all $m\in\lC-r,r-1\rC\cap\Z$. 
Write $I^r_m$ and $J^r_{m}$ the corresponding arcs on the unit circle such that their projections on the horizontal axis is the interval $\lC \frac{m}{r}, \frac{m+1}{r}\rC$. 
If $4r\leq n$, then the Pigeon-hole principle implies that there exists at least one $M\in\lC-r,r-1\rC\cap\Z$ such that $I_{M}$ or $J_M$ contains at least $\frac{n}{4r}\geq 1$ elements of the set  $\lL \exp\lP i \lP j 2\pi x_\alpha - \theta \rP\rP : j=0,\ldots, n-1\rL$.

In the sequel, we define
\[
I^r_M:=\left\{j\in \{0,\ldots,n-1\}:\cos\lP j2\pi x - \theta\rP\in\lC\frac{M}{r},\frac{M+1}{r}\rC\in 
\lC\frac{M}{r},\frac{M+1}{r}\rC\right\}\not=\emptyset.
\]
and for each $j\in I^r_M$, we define
\[
d_j = \min\lL \left|\cos\lP j 2\pi x_\alpha - \theta\rP - \frac{M}{r}\right|, \left|\cos\lP j 2\pi x_\alpha - \theta\rP - \frac{M+1}{r}\right|\rL.\]
Note that 
\[
\min_{0\;\leq\; l<k\; \leq\; n-1} \abs{l 2\pi x_\alpha - k 2\pi x_\alpha} \geq  \frac{2\pi\alpha'}{N'}\geq  \frac{2\pi}{N'}.
\]
Let $L =\min\lL \left\lfloor \frac{n}{8r} - \frac{3}{2} \right\rfloor, \left\lfloor \frac{N' }{8r} -\frac{1}{2}\right\rfloor \rL$ and observe that 
for each $0\leq\lambda\leq L$ there exists at least $j\in I^r_M$ such that $d_j\geq \lP 2\lambda + 1\rP  \frac{2\pi}{N'}$. Then 
\[
\begin{split}
s^r_M:=\sum_{j\in I^r_M}d_j\geq \sum_{\lambda=0}^{L} \lP 2\lambda + 1\rP  \frac{2\pi}{N'}=\frac{2\pi(L+1)^2}{N'}\geq \frac{2\pi L^2}{N'}.
\end{split}
\]
By the choosing of $L$, if $r \leq \left\lfloor n^{1/4} \right\rfloor$, we have
$
\frac{2\pi L^2}{N'} \geq \frac{ 2\pi}{n^{11/10}}n^{3/2}
$
for all large $n$.

Here, let $v$ be a vector in $\R^n$ with entries  $v_j=\cos\lP j 2\pi x_\alpha - \theta\rP$ for each $j=0,\ldots,n-1$. If a positive integer $r\leq \left\lfloor n^{1/4} \right\rfloor$, then by the previous discussion we deduce 
$\mathrm{dist}(rv,\Z^n)\geq  2\pi n^{2/5}$ for all $n$ large. If $r$ is any positive real number, observe that $\left[ \frac{s}{r},\frac{s+1}{r} \right]\subset \left[ \frac{s}{\lceil r\rceil},\frac{s+1}{\lceil r\rceil} \right]$,
where $s\in \mathbb{N}$, and therefore our previous analysis holds true for any $r>0$.

Assume that $r\leq \left\lfloor n^{1/4}\right\rfloor$ and recall that $\|V_\alpha^T \Theta\|_2\leq r\sqrt{n}$ and that
$L\geq \sqrt{2/q}$ is fixed.
By the definition of lcd, for all $n$ large we obtain
\[
\begin{split}
2\pi n^{2/5} \leq   \textnormal{dist}\lP V_\alpha^T \Theta, \Z^n\rP 
\leq  L\sqrt{ \log_+ \frac{\|{V_\alpha^T \Theta}_2\|_2}{L} } 
 \leq L\sqrt{\log_+\frac{n^{3/4}}{L}},
\end{split}
\]
which yields a contradiction for $n$ large. Thus, 
$D\lP V_\alpha\rP \geq \left\lfloor n^{1/4}\right\rfloor$. 
Therefore, the preceding inequality together with inequality \eqref{eqn22022019901}, inequality   
\eqref{eqn110320191612} and the fact the cardinality of 
$J_3(n,N)$ is at most $n^{11/10+\textnormal{o}(1)}(\log(n))^{\textnormal{o}(1)}$ 
allow us to deduce
\begin{eqnarray*}
& &  \sum\limits_{\alpha \in J_3(n,N)}   \Prob{\abs{W_n\lP e^{i2\pi x_\alpha} \rP} \leq 2tC_2 n^{-1/2}\lP\log n\rP^{-2}}  \\ [0.3em]
& &\leq n^{11/10+\textnormal{o}(1)}(\log(n))^{\textnormal{o}(1)} \lP \frac{2C^2L^2\lP 2tC_2\rP^2}{\lP c_4n^{2-1/5}\rP^{1/2}\lP n^{1/2}\lP\log n\rP^2\rP^2}\rP \\ [0.3em]
& & \;\;\;\;+\;n^{11/10+\textnormal{o}(1)}(\log(n))^{\textnormal{o}(1)} \lP\frac{2C^2L^2}{\lP c_4n^{2-1/5}\rP^{1/2} \lP n^{1/4} \rP^2} \rP  \\ [0.3em]
& & \leq \;\frac{8C^2C^2_2L^2 t^2 }{c_4^{1/2} n^{4/10} }+ \frac{2 C^2 L^2 }{c_4^{1/2} n^{1/10}}
\end{eqnarray*}
for all large $n$.
As a consequence we obtain
\begin{equation*}
\sum\limits_{\alpha \in J_3(n,N)}  \Prob{\abs{W_n\lP e^{i2\pi x_\alpha} \rP} \leq 2tC_2 n^{-1/2}\lP\log n\rP^{-2}} = \textnormal{O}\lP\frac{1}{n^{1/10}}\rP, 
\end{equation*} where the implicit constant depends on $L$ and $t$.
\end{proof}

{
\appendix
\label{appendix}
\section{Arithmetic properties}\label{A:app270120191754}
This section contains the proofs of the results that we skipped in the paper in order to be more fluid.

\begin{lemma}\label{app270120191754}
Let $m\geq 1$ and $M\in \mathbb{N}$. Then 
the cardinality of the set
\[
\Gamma^M_m:=\left\{k\in [1,M]\cap \mathbb{N}:~\gcd\lP k, M\rP \geq m\right\}
\] 
is at most $\frac{1}{\lfloor m\rfloor} M^{1+C\lP \log \log M\rP^{-1}}$, where $C$ is a positive constant.
\end{lemma}
\begin{proof}
Denote by $T$ the Euler totient function. Observe that 
\[
\sum\limits_{k\in \Gamma^{M}_m}   1 
\leq 
\sum_{\substack{
        d=\lfloor m\rfloor\\
        d\left| M\right.
    } 
}^M   T\lP\frac{M}{d}\rP. 
\]
It is well-known that  $T\lP s\rP\leq s - \sqrt{s}$ for all $s\in\N$. Moreover, if $d(s)$ denotes the number of positive divisors of $s$, then  Theorem 13.12 in \cite{APO} implies that  there exists a positive  constant $C$ such that $d(s) \leq s^{C\lP\log\log\lP s\rP\rP^{-1}}$.
Hence,
\begin{eqnarray} \label{eqn220220191514}
\sum_{k\in \Gamma^{M}_m} 1 
 \leq 
\lP\frac{M}{\lfloor m\rfloor} - \sqrt{\frac{M}{\lfloor m\rfloor}}\rP M^{C\lP\log\log\lP M\rP\rP^{-1}}  \nonumber 
 \leq  \frac{1}{\lfloor m\rfloor} M^{1+C\lP \log \log M\rP^{-1}}
\end{eqnarray}
which yields the statement. 
\end{proof}

\begin{lemma}\label{lowerbound}
Let $\theta\in[0,2\pi)$ and $n\in\mathbb{N}$.
Let $\mathcal{V}=(\mathcal{V}_j)_{j\in \{1,\ldots,n\}}\in \R^n$ such that $\mathcal{V}_j= r\cos\lP j 2\pi x-\theta \rP$ for $j=0,\ldots,n-1$, where  $r\in \N$ and $x={1}/{n}$. 
If $\frac{1}{4\pi r x}\geq 8$, then
\[
\mathrm{dist}\lP \mathcal{V},\Z^n\rP \geq \frac{1}{128\pi x }.
\]
\end{lemma}
\begin{proof}
Let $\theta\in[0,2\pi)$ and $n\in\mathbb{N}$. Let $x={1}/{n}$ and we define the following sequence $P_n=\lL\exp\lP i\lP j2\pi x - \theta\rP\rP: j=0,\ldots,n-1\rL$, where $i$ is imaginary unity. Note $P_n$ is a set of points on the unit circle which can be looked as vertices of a regular polygon with $n$ sides inscribed in the unit circle. 

Since the arguments of ``two consecutive points" on $P_n$, $\exp\lP i\lP j2\pi x - \theta\rP\rP$ and $\exp\lP i\lP (j+1)2\pi x - \theta\rP\rP$, are separated by a distance $2\pi x$, the number of points in $P_n$ which are in any arc of length $\ell $ on the unit circle is at least $\frac{\ell}{2\pi x}-2$.

Let $\lC y, y + 8\pi x\rC$ be a subinterval of $[-1,1]$.  We consider an arc $\stackrel{\frown}{I}$
 on the unit circle such that its projection on the horizontal axis is $\lC y, y + 8\pi x\rC$. If the length of the arc 
 $\stackrel{\frown}{I}$
 is $\ell$, then the number of values $\cos\lP j2\pi x - \theta\rP$, $j=0,\ldots,n-1$ that belongs to $\lP y, y + 8\pi x\rP$ is at least $\frac{1}{2}\left(\frac{\ell}{2\pi x}-2\right)$. 
Observe that  $\frac{1}{2}\left(\frac{\ell}{2\pi x}-2\right)\geq 1$ when 
$\ell\geq 8\pi x$. 

Let $r\in \mathbb{N}$ and $m\in\lC-(r-1),(r-1)\rC\cap\Z$. By the preceding explanation, 
for all positive integer $k\leq \frac{1}{8\pi rx}$
there exists $j\in \{0,\ldots,n-1\}$
such that
\[
\cos\lP j2\pi x - \theta\rP\in\lP \frac{m}{r} + 8\pi\lP k-1\rP
 x ,\frac{m}{r} + 8\pi k x\rP \subset \lC\frac{m}{r},\frac{m+1}{r}\rC.
\] 
In the sequel, set
\[
I^r_m:=\left\{j\in \{0,\ldots,n-1\}:\cos\lP j2\pi x - \theta\rP\in 
\lC\frac{m}{r},\frac{m+1}{r}\rC\right\}\not=\emptyset
\]
and for each $j\in I^r_m$, define
\[d_j:=\min\lL \left|\cos\lP j2\pi x - \theta\rP - \frac{m}{r}\right|, \left|\cos\lP j2\pi x - \theta\rP - \frac{m+1}{r}\right|\rL.
\] 
Let $L$ be the biggest integer such that $8\pi Lx \leq \frac{1}{2r}$, or equivalently, $L=\left\lfloor \frac{1}{16 \pi rx} \right\rfloor$. 
Observe that
\[
L \geq \frac{1}{16\pi rx} - 1 \geq \frac{1}{32\pi rx}\quad \textrm{ when }\quad \frac{1}{4\pi rx}\geq 8.
\]
Then 
\[
s^r_m:=\sum\limits_{j\in I^r_m} d_j\geq \sum\limits_{\lambda=1}^L 2\lambda (8\pi x)\geq 8\pi x L^2\geq \frac{1}{128\pi r^2x}.
\]
Moreover,
\begin{equation*}\label{eqn280221091800}
\sum_{m=-(r-1)}^{m=r-1} s^r_m \geq  \frac{2r-1}{128\pi r^2x} \geq \frac{1}{128\pi rx}, 
\end{equation*} 
where the last inequality follows since $\frac{2r-1}{r}\geq 1$ for $r\in \mathbb{N}$.
Consequently,  the distance between the vector $\mathcal{V}\in\R^n$ with  entries $\mathcal{V}_j = r\cos\lP j2\pi x - \theta\rP$ for $j=0,\ldots,n-1$ with $x={1}/{n}$ and the set $\Z^n$ is at least
\begin{equation*}\label{eqn280120191810}
r\left(\frac{1}{128\pi rx}\right)=\frac{1}{128\pi x}\quad \textrm{ verifying that } \frac{1}{4\pi rx} \geq 8 \textrm{ is fulfilled}.
\end{equation*} 
\end{proof}
}
\section{Acknowledgments}
Both authors are grateful with professor Jes\'us L\'opez Estrada for his reading and recommendations in the preliminary version of this work.
G. Barrera acknowledges support from a post-doctorate grant held at
Center for Research in Mathematics, (CIMAT, 2015--2016). He would like to express his gratitude to Pacific Ins\-titute for the Mathematical Sciences (PIMS, 2017--2019) for the grant held
at the Department of Mathematical and Statistical Sciences at University of Alberta. He also would like to thank CIMAT and University of Alberta for all the facilities used along the realization of this paper.
P.~Manrique acknowledges support 
from C\'atedras CONACYT for the 
research position held at
Mathematics Institute, Cuernavaca (UNAM, 2017--2020).
He also would like to thank UNAM for all the facilities used along the realization of this paper.
Both authors would like to thank the constructive and useful suggestions provided by the 
anonymous referee.
\markboth{}{References}
\bibliographystyle{amsplain}

\end{document}